\documentclass[1p]{elsarticle}
\usepackage{amssymb}
\usepackage{amsmath,amssymb,amstext}
\usepackage{amsfonts}
\usepackage{float}
\usepackage[figuresright]{rotating}
\usepackage[hidelinks]{hyperref}
\usepackage[dvipsnames]{xcolor}
\usepackage{tikz}
\usepackage{algorithm}
\usepackage{algorithmicx}
\usepackage{algpseudocode}
\usepackage{pifont}
\usepackage{mathtools}
\usepackage{float}

\usepackage{array}
\usepackage{eqparbox}
\usepackage{fancyvrb}
\usepackage{graphicx}

\usepackage{lscape}
\usepackage{newlfont}

\usepackage{caption}
\usepackage{subcaption}

\usepackage{amsthm}
\newtheorem{theorem}{Theorem}[section]
\newtheorem{lemma}{Lemma}[section]

\newtheorem{definition}{Definition}[section]

\newtheorem{remark}{Remark}[section]

\newcommand{\blue}[1]{#1}

\journal{Numerical Methods for PDEs}

\begin{document}

\begin{frontmatter}
\title{Fourier Spectral Methods with Exponential Time Differencing for Space-Fractional Partial Differential Equations in Population Dynamics} 
	\date{July 22, 2022} 
\author[]{A. P. Harris\corref{mycorrespondingauthor}}
\cortext[mycorrespondingauthor]{Corresponding author}
\ead{ashlin_harris@brown.edu}
\address{\blue{Brown Center for Biomedical Informatics},
 Brown University, Providence, RI, USA}
\author{T. A. Biala}
\address{Department of Mathematics, The Ohio State University, Columbus, OH, USA}
\author{A. Q. M.  Khaliq}
\address{Department of Mathematical Sciences and Center for Computational Science,
 Middle Tennessee State University, Murfreesboro, TN, USA}

\begin{abstract}
Physical laws governing population dynamics are generally expressed as differential equations.
Research in recent decades has incorporated fractional-order (non-integer) derivatives into differential models of natural phenomena, such as
reaction-diffusion systems.

In this paper, we develop a method to numerically solve a multi-component and multi-dimensional space-fractional system.
For space discretization, we apply a Fourier spectral method that is suited for multidimensional PDE systems.
Efficient approximation of time-stepping is accomplished with a locally one dimensional exponential time differencing approach.
We show the effect of different fractional parameters on growth models
and consider the convergence, stability, and uniqueness of solutions,
as well as the biological interpretation of parameters and boundary conditions.

\end{abstract}
\end{frontmatter}

\section{Introduction}

\blue{Diffusion, the gradient-driven net movement of a substance, occurs at widely diverse scales in nature.
It is a topic of study in all the natural sciences and has been modelled in terms of both differential equations and random walk simulations.}
Reaction-diffusion (RD) models have been developed to
describe the movement and interaction of species in various contexts, especially chemistry and biology.
They have been used extensively in the study of ecological predator-prey
interactions
and have been called reproduction-dispersal models in this context\cite{ashi2018solving}.

\autoref{equation:classical RD} shows the classical reaction-diffusion equation,
\blue{where $\vec{D}$ is a diagonal matrix of diffusion components.
Population density $u$ is given over the spatial component $x$ and temporal component $t$:}
\begin{equation}
\label{equation:classical RD}
\dfrac{\partial u}{\partial t} = \vec{D}\dfrac{\partial ^2 u}{\partial x^2} + f(u).
\end{equation}
The function $f$, the reaction term, stands for population growth; Fisher's equation, defined as $f(u)=ru(1-u/K)$, is widely used\cite{baeumer2007fractional}.
For that choice of reaction term, $r$ stands for the intrinsic growth rate of the species, and $K$ is the carrying capacity.

In practice, biologists may introduce other parameters and functions that capture the complexities of \blue{interactions}.
In applications of reaction-diffusion models to predator-prey interactions, \blue{the system is influenced primarily by the functional response,  the rate at which an individual predator consumes prey.
The diffusion terms of nonlinear reaction diffusion problems tend to be stiff,
constraining the step size of any fully explicit numerical solution\cite{ashi2018solving}.}

Recently, attention has been given to anomalous diffusion in nature, which features a stable nonlinear relationship between mean squared displacement and time.
\blue{In the context of physics and chemistry, anomalous diffusion is generally described as a fractional-order movement across space or time.
Alternatively, biologists might describe anomalous diffusion in terms of a fractional-order dimension of the time or space;
that is, the entity exhibits classical diffusion over a fractal space\cite{johnson1992diffusion}.
For instance, Baeumer \textit{et al.} have proposed the use of fractional-in-space models to capture realistic spreading behaviour of invading species\cite{baeumer2007fractional}.}

\blue{\autoref{equation:classical RD} may be generalised to include a fractional (non-integer) differential operator.}
The new equation, as defined in \cite{baeumer2007fractional},
is given by
\begin{equation}
\label{equation:fractional RD}
\dfrac{\partial u}{\partial t} = \vec{D}\dfrac{\partial ^\alpha u}{\partial x^\alpha} + f(u),
\end{equation}
where the value of $\alpha$ is constrained to $1 < \alpha \leq 2$.
The special case $\alpha=2$ yields the classical reaction-diffusion equation, as in \autoref{equation:classical RD}.

Anomalous transport has been observed in complex polymer networks and materials with varying sizes of pores or obstacles.
One characteristic of systems with anomalous diffusion is macromolecular crowding\cite{hofling2013anomalous}.
The cytoplasm of cells is densely packed with molecules of vastly different sizes and interactions.
Similarly, in materials with various pore sizes, diffusing elements must react with a complex substrate surface.

Interestingly, the cell interior gives examples of anomalous diffusion of molecules in
1-D (the transport of molecules and vesicles along microfilaments),
2-D (diffusion across plasma membranes),
and 3-D space (diffusion of molecules in crowded cytoplasm)\cite{krapf2015mechanisms,regner2013anomalous,regner2014identifying}.
Instances of cellular diffusion include
cell migration,
neuron growth and neural crest development,
and metastasis.

Anomalous diffusion can be modelled in terms of fractional-order differential equations and L\'{e}vy flights.
However, numerical methods developed for integer-order systems are not typically suited for the efficient computation of fractional-order systems, since they feature dense matrix operations that become expensive for systems with three or more dimensions.
Additionally, differential equation models of biological systems can include many components subject to nonlinearities and memory effects, requiring efficient numerical methods that scale suitably.

\section{Methods}

\begin{definition}
\label{def:Laplacian}
Suppose the Laplacian $(- \Delta)$ has a complete set of orthonormal eigenfunctions $\phi_n, ~\phi_{n,m},~ \text{or}~ \phi_{n,m,l}$ corresponding to the eigenvalues $\lambda_n, ~\lambda_{n,m},~\text{or}~\lambda_{n,m,l},$ respectively, on a bounded region $\Omega$. \blue{That is, for $n,m,l = 0,1,2, \cdots$ in $\Omega$},
\[(-\Delta) \phi_n = \lambda_n\phi_n , ~~~~~~~~~~~~~~~d=1,\]
\[(-\Delta) \phi_{n,m} = \lambda_{n,m}\phi_{n,m} , ~~~~~~d=2,\]
\[ (-\Delta) \phi_{n,m,l} = \lambda_{n,m,l}\phi_{n,m,l}, ~~ d=3,\]
and ${\mathcal{B}}(\phi) = 0$ on $\partial \Omega$, where $\mathcal{B}(\phi)$ are the homogenous Dirichlet or homogenous Neumann boundary conditions. Let 
\begin{equation*}
\begin{split}
f_1 &= \sum_{n=0}^{\infty}c_n\phi_n~\text{ such that }~ \sum_{n=0}^{\infty}|c_n|^2|\lambda_{n}|^{\alpha} < \infty, ~~~~~~~~~~~~~~~~~~~~~~~~~~~~~~~~~~~~~~~~~d=1,\\
f_2 &= \sum_{n=0}^{\infty}\sum_{m=0}^{\infty}c_{n,m}\phi_{n,m}~\text{ such that }~ \sum_{n=0}^{\infty}\sum_{m=0}^{\infty}|c_{n,m}|^2|\lambda_{n,m}|^{\alpha} < \infty, ~~~~~~~~~~~~~~~~~~~~~d=2,\\
f_3 &= \sum_{n=0}^{\infty}\sum_{m=0}^{\infty}\sum_{l=0}^{\infty}c_{n,m,l}\phi_{n,m,l}~\text{ such that }~ \sum_{n=0}^{\infty}\sum_{m=0}^{\infty}\sum_{l=0}^{\infty}|c_{n,m,l}|^2|\lambda_{n,m,l}|^{\alpha} < \infty, ~~~~~~~d=3.
\end{split}
\end{equation*}
 Then, 
 $(-\Delta)^{\frac{\alpha}{2}}$ is defined by
\begin{equation*}
\begin{split}
(-\Delta)^\frac{\alpha}{2}f_1 &= \sum_{n=0}^{\infty}c_n\lambda_{n}^{\frac{\alpha}{2}}\phi_n, ~~~~~~~~~~~~~~~~~~~~~~~~~~~~~~~~~~~~~~~~~d=1,\\
(-\Delta)^\frac{\alpha}{2}f_2 &= \sum_{n=0}^{\infty}\sum_{m=0}^{\infty}c_{n,m}\lambda_{n,m}^{\frac{\alpha}{2}}\phi_{n,m}, ~~~~~~~~~~~~~~~~~~~~~~~~~~~~d=2,\\
(-\Delta)^\frac{\alpha}{2}f_3 &= \sum_{n=0}^{\infty}\sum_{m=0}^{\infty}\sum_{l=0}^{\infty}c_{n,m,l}\lambda_{n,m,l}^{\frac{\alpha}{2}}\phi_{n,m,l}, ~~~~~~~~~~~~~~~~~~~~d=3.
\end{split}
\end{equation*}
\end{definition}
\begin{remark}
 For homogeneous Dirichlet boundary conditions with 

$\Omega = (a, b)^d, ~d=1,2,3,$ and $\mathbf{x} \in \Omega$,
 
\[\lambda_{\eta_1, \cdots, \eta_d} = \sum_{n=\eta_1}^{\eta_d}\left(\dfrac{(n+1)\pi}{b-a}\right)^2,\]

 \[\phi_{\eta_1, \cdots, \eta_d} = \left(\sqrt{\dfrac{2}{b-a}}\right)^d\prod_{n=\eta_1}^{\eta_d}\sin\left(\dfrac{(n+1)\pi (x_n - a)}{b-a}\right),~~~\eta_i = 0,1,2, \cdots.\]
 For homogeneous Neumann boundary conditions, 
\[\lambda_{\eta_1, \cdots, \eta_d} = \sum_{n=\eta_1}^{\eta_d}\left(\dfrac{n\pi}{b-a}\right)^2,\]
and 
\[\phi_{\eta_1, \cdots, \eta_d} = \left(\sqrt{\dfrac{2}{b-a}}\right)^d\prod_{n=\eta_1}^{\eta_d}\cos\left(\dfrac{n\pi (x_n - a)}{b-a}\right), ~~~\eta_i = 0,1,2, \cdots.\]
\end{remark}
\subsection{Spatial Discretization}
We begin with the reaction-diffusion system of $M$ species $\mathbf{u} = [u_1, u_2, \cdots, u_M]^T$,
\begin{equation}
    \label{eqn:system1}
    \dfrac{\partial{u}_i}{\partial t} = -D_i\left( -\Delta\right)^{\alpha/2}u_i + f_i(\textbf{u}), ~~~~~~ 1 \leq \alpha \leq 2, ~~~~\textbf{x} \in \mathbb{R}^n,
\end{equation}
\blue{where $D_i$ is the diffusion coefficient of the $i$\textsuperscript{th} specie $u_i$} and $f_i(\textbf{u})$ is its reaction term.
By applying a Fourier transform and the definition of the fractional Laplacian (\ref{def:Laplacian}) to (\ref{eqn:system1}), we obtain (for $d=1$)
\begin{equation}
    \label{eqn:system2}
    \dfrac{\partial{\hat{u}}_{ij}}{\partial t} = -D_i\lambda_j^{\alpha/2}\hat{u}_{ij} + \hat{f}_{ij}(\hat{\textbf{u}}), 
\end{equation}
where $\hat{u}_{ij}$ is the \blue{$j$\textsuperscript{th}} Fourier coefficients of the \blue{$i$\textsuperscript{th}} specie and $\hat{f}_{ij}$ is its associated reaction term.
Note that \blue{the orthogonality of the basis functions implies that} each of the Fourier coefficients evolve independently of one another.
\blue{For $L$ grid points} in space ($j = 1, \cdots, L$) and step size $h$, the corresponding homogeneous Dirichlet boundary conditions are
\[ \lambda_j = \dfrac{j\pi}{b-a}, ~~~ x_j = a + jh + \dfrac{h}{2}, ~~ h = \dfrac{(b-a)}{L+1},\]
and the homogeneous Neumann boundary conditions are given by
\[ \lambda_j = \dfrac{(j-1)\pi}{b-a}, ~~~ x_j = a + (j - 1)h + \dfrac{h}{2}, ~~ h = \dfrac{(b-a)}{L}.\]
For more information about the mesh generation, we refer the reader to \cite{bueno2014fourier}. We compute the coefficients $\hat{u}_i$ and the inverse reconstruction of $u$ \blue{in physical space using} coefficient algorithms (discrete sine or cosine transforms and their inverses)
based on the specified homogeneous boundary conditions\cite{bueno2014fourier,briggs1995dft, bueno2006spectral}.
\subsection{Time Discretization}
\blue{To discretize across time, }we rewrite the system (\ref{eqn:system2}) as 
\begin{equation}
    \label{eqn:system3}
    \dfrac{\partial{\hat{U}}_{i}}{\partial t} = -D_i\Lambda^{\alpha/2}\hat{U}_{i} + \hat{F}_{i}(\hat{\textbf{U}}),
\end{equation}
where \blue{$\hat{U}_i$ and $\hat{F}_i$
 are the Fourier coefficient and reaction term respectively for the $i$\textsuperscript{th} specie.}
Let $t_k=k\tau, ~k= 0,...,N$, where $\tau = T/N$ is the time step size, and $\hat{U}(t_k) := \hat{U}^k$.
Here, the exact solution of (\ref{eqn:system3}) at time $t_{k+1}$ can be written as 
\begin{equation}
\label{eqn:generalETDIntegral}
\hat{U}_i(t_{k+1}) = e^{-\tau {\Lambda^{\frac{\alpha}{2}}}}\hat{U}_i(t_k) + \tau\int_{0}^{1}e^{-\tau  {\Lambda^{\frac{\alpha}{2}}}\left(1-s\right)}\hat{F}_i(\hat{\mathbf{U}}(t_k+s\tau))ds,
\end{equation}
\blue{where $\hat{\mathbf{U}}(t_k) = [\hat{U}_1(t_k), \hat{U}_2(t_k), \cdots, \hat{U}_M(t_k)]^T$}. Equation (\ref{eqn:generalETDIntegral}) serves as the basis for exponential time differencing (ETD) schemes which are obtained by using different approximations to the matrix exponential function and the nonlinear reaction terms. Suppose \blue{$\hat{F}_i(\hat{\mathbf{U}}(t_k + s\tau))$,} in  (\ref{eqn:generalETDIntegral}), is approximated by an average over end points in an interval $[t_k, t_{k+1}]$. That is, 
\[\hat{F}_i(\hat{\mathbf{U}}) \approx \hat{F}_i(\hat{\mathbf{U}}^k) + (t - t_k)\dfrac{\hat{F}_i(\hat{\mathbf{b}}^k)- \hat{F}_i(\hat{\mathbf{U}}^{k})}{\tau}, ~~~~ t \in[t_k, t_{k+1}],\]
where 
\begin{align*}
\blue{\hat{\mathbf{b}}^k }&= \blue{[\hat{b}_1^k, \hat{b}_2^k, \cdots, \hat{b}_M^k]^T,}\\
\hat{b}^k_i &= e^{-\tau  \Lambda^{\frac{\alpha}{2}}}\hat{U}_i^k  + \Lambda^{-\frac{\alpha}{2}}\left({I} - e^{-\tau  \Lambda^{\frac{\alpha}{2}}}\right)\hat{F}_i(\hat{\mathbf{U}}^{k}).
\end{align*}
Then, the integral equation (\ref{eqn:generalETDIntegral}) becomes
\begin{equation*}
\hat{U}_i^{k+1} \approx e^{-\tau  \Lambda^{\frac{\alpha}{2}}}\hat{U}_i^k + \tau e^{-\tau  \Lambda^{\frac{\alpha}{2}}}\int_{0}^{1}e^{\tau  \Lambda^{\frac{\alpha}{2}}s}\left(\hat{F}_i(\hat{\mathbf{U}}^{k}) + s\tau\dfrac{\hat{F}_i(\hat{\mathbf{b}}^{k})- \hat{F}_i(\hat{\mathbf{U}}^{k})}{\tau}\right)ds,
\end{equation*}
which simplifies to
\begin{equation}
\label{eqn:ETDCNtemp}
\hat{U}_i^{k+1} = \hat{b}_i^k + \dfrac{1}{\tau}\Lambda^{-\alpha} \left(e^{-\tau  \Lambda^{\frac{\alpha}{2}}}-I+\tau  \Lambda^{\frac{\alpha}{2}}\right)\left[\hat{F}_i(\hat{\mathbf{b}}^{k})  - \hat{F}_i(\hat{\mathbf{U}}^{k})\right],
\end{equation}
\blue{where $I = [1, \cdots, 1]^T$}. Replacing the exponential matrix in  (\ref{eqn:ETDCNtemp}) by the (1,1)-Pad\'e approximation, the Crank-Nicolson ETD (ETD-CN) method is obtained, after some simplification, as
\begin{equation}
\begin{split}
\label{eqn:ETDCN}
\hat{V}_i^{k+1} &= \hat{a}_i^k + \tau \left(2{I} + \tau  \Lambda^{\frac{\alpha}{2}}\right)^{-1}\left[\hat{F}_i(\hat{\mathbf{a}}^{k}) - \hat{F}_i(\hat{\mathbf{V}}^{k})\right]\approx \hat{U}_i^{k+1},\\
\hat{a}_i^k &= \left\{4\left(2{I} + \tau  \Lambda^{\frac{\alpha}{2}}\right)^{-1}-{I}\right\}\hat{V}_i^k + 2\tau \left(2{I} + \tau  \Lambda^{\frac{\alpha}{2}}\right)^{-1}\hat{F}_i(\hat{\mathbf{V}}^{k}),
\end{split}
\end{equation}
where $\hat{\mathbf{a}}^k = [\hat{a}_1^k, \hat{a}_2^k, \cdots, \hat{a}_M^k]^T$ and $\hat{\mathbf{V}}^k = [\hat{V}_1^k, \hat{V}_2^k, \cdots, \hat{V}_M^k]^T$. 
 \blue{\subsection{Error Analysis}
 Here, we discuss the convergence and stability of the numerical scheme (\ref{eqn:ETDCN}). We assume that the nonlinear function $\hat{\mathbf{F}}(\hat{\mathbf{U}})$ is Lipschitz continuous in a region   $\Omega \times (0,T]~ (\Omega \subset \mathbb{R}^n)$.
 That is, there exists a constant $L$ such that for $\hat{\mathbf{U}},~\hat{\mathbf{V}}~\text{in}~ \Omega \times (0,T]$, 
\[||\hat{\mathbf{F}}(\hat{\mathbf{U}}) - \hat{\mathbf{F}}(\hat{\mathbf{V}})|| \leq L||\hat{\mathbf{U}} - \hat{\mathbf{V}}||,\]
where $||\cdot||$ is the $\ell_2$-norm. Lipschitz continuity on $\hat{\mathbf{F}}(\hat{\mathbf{U}})$  implies that
\begin{equation}
\label{eqn:LipschitzAssumption}
||\hat{\mathbf{F}}(\hat{\mathbf{U}})||  -  ||\hat{\mathbf{F}}(\hat{\mathbf{V}})||  \leq ||\hat{\mathbf{F}}(\hat{\mathbf{U}}) - \hat{\mathbf{F}}(\hat{\mathbf{V}})||  \leq L||\hat{\mathbf{U}} - \hat{\mathbf{V}}||
\end{equation}
and
\begin{equation}
\label{eqn:LipschitzAssumption2}
||\hat{\mathbf{F}}(\hat{\mathbf{U}})||  \leq L||\hat{\mathbf{U}}|| + D,
\end{equation}
where $D = ||\hat{\mathbf{F}}\mathbf{(0)}||$. In the following analysis, we take $K$ to be a generic positive constant. \\
\begin{lemma}
\label{lem:UnEstimate2}
Suppose $\hat{\mathbf{F}}(\hat{\mathbf{U}})$ is Lipschitz continuous. Then,
\[||\hat{\mathbf{b}}^k - \hat{\mathbf{U}}^k|| \leq  K\tau\left( ||\hat{\mathbf{U}}^k|| + D\right).\]
\end{lemma}
\begin{proof}
\begin{equation*}
\begin{split}
||\hat{\mathbf{b}}^k - \hat{\mathbf{U}}^k|| &\leq ||e^{-\tau \Lambda^{\frac{\alpha}{2}}} - I||\cdot||\hat{\mathbf{U}}^k|| + \left|\left|\Lambda^{-\frac{\alpha}{2}}\left(e^{-\tau \Lambda^{\frac{\alpha}{2}}} - I\right)\right|\right|\cdot||\hat{\mathbf{F}}(\hat{\mathbf{U}}^k)||\\
&\leq K\tau||\hat{\mathbf{U}}^k|| + K\tau\left(L||\hat{\mathbf{U}}^k|| +  D\right),
\end{split}
\end{equation*}
and so the result must follow for some constant $L$.
\end{proof}
\begin{lemma}
If $\hat{\mathbf{F}}(\hat{\mathbf{U}})$ is Lipschitz continuous, then the sequence of solutions defined by (\ref{eqn:ETDCNtemp}) satisfies
\[||\hat{\mathbf{U}}^k|| \leq K\left(||\hat{\mathbf{U}}^0|| + D\right),\]
where $D$ is as defined in (\ref{eqn:LipschitzAssumption2}).
\end{lemma}
\begin{proof}
A recursive application of (\ref{eqn:ETDCN}) gives
\begin{equation*}
\begin{split}
\hat{\mathbf{b}}^{k-1} & = e^{-k\tau  \Lambda^{\frac{\alpha}{2}}}\hat{\mathbf{U}}^{0} + \dfrac{1}{\tau}\Lambda^{-\alpha}\left(e^{-\tau  \Lambda^{\frac{\alpha}{2}}} - I +  \tau \Lambda^{\frac{\alpha}{2}}\right) \sum_{j=1}^{k-1}e^{-j\tau \Lambda^{\frac{\alpha}{2}}}\left[\hat{\mathbf{F}}(\hat{\mathbf{b}}^{k-j-1}) -  \hat{\mathbf{F}}(\hat{\mathbf{U}}^{k-j-1})\right]\\
&~~ +  \Lambda^{-\frac{\alpha}{2}}\left(I - e^{-\tau  \Lambda^{\frac{\alpha}{2}}}\right)\sum_{j=0}^{k-1}e^{-j\tau \Lambda^{\frac{\alpha}{2}}}\hat{\mathbf{F}}(\hat{\mathbf{U}}^{k-j-1}),
\end{split}
\end{equation*}
which implies that
\begin{equation*}
||\hat{\mathbf{b}}^{k-1}|| \leq K||\hat{\mathbf{U}}^{0}|| + K\tau\sum_{j=1}^{k-1}||\hat{\mathbf{b}}^{k-j-1} - \hat{\mathbf{U}}^{k-j-1}|| + K\tau\left(\sum_{j=0}^{k-1}||\hat{\mathbf{U}}^{j}|| + D\right).
\end{equation*}
Now,
\begin{equation*}
\begin{split}
||\hat{\mathbf{U}}^{k}|| &\leq ||\hat{\mathbf{b}}^{k-1}|| + \dfrac{1}{\tau}||\Lambda^{-\alpha}\left(e^{-\tau \Lambda^{\frac{\alpha}{2}}} - I + \tau \Lambda^{\frac{\alpha}{2}}\right)||\cdot||\hat{\mathbf{F}}(\hat{\mathbf{b}}^{k-1})- \hat{\mathbf{F}}(\hat{\mathbf{U}}^{k-1})||\\
& \leq K||\hat{\mathbf{U}}^{0}|| + K\tau\sum_{j=1}^{k-1}||\hat{\mathbf{b}}^{k-j-1} - \hat{\mathbf{U}}^{k-j-1}|| + K\tau\left(\sum_{j=0}^{k-1}||\hat{\mathbf{U}}^{j}|| +  D\right) + K\tau~||\hat{\mathbf{b}}^{k-1} - \hat{\mathbf{U}}^{k-1}||\\
& = K||\hat{\mathbf{U}}^{0}|| + K\tau\sum_{j=0}^{k-1}||\hat{\mathbf{b}}^j - \hat{\mathbf{U}}^{j}|| + K\tau\left(\sum_{j=0}^{k-1}||\hat{\mathbf{U}}^{j}|| + D\right).
\end{split}
\end{equation*}
Using the error bound in Lemma (\ref{lem:UnEstimate2}) we obtain
\begin{equation*}
\begin{split}
||\hat{\mathbf{U}}^k|| & \leq K||\hat{\mathbf{U}}^{0}|| + K\tau\left(\sum_{j=0}^{k-1}||\hat{\mathbf{U}}^{j}|| + D\right) + K\tau^2\left(\sum_{j=0}^{k-1}||\hat{\mathbf{U}}^{j}|| + D\right).
\end{split}
\end{equation*}
The rest of the proof is established by Gronwall's lemma for the discrete case\cite[pp. 8]{Emmrich1999}.
\end{proof}
\begin{lemma}
\label{lem:UnEstimate3}
The following error bound 
\[||\hat{\mathbf{b}}^{k-1} - \hat{\mathbf{a}}^{k-1}|| \leq   K\tau^2||\hat{\mathbf{U}}^{0}||  + K\tau^3\sum_{j=0}^{k-1}||\hat{\mathbf{U}}^{j}||+ K\tau^3D + K\tau\sum_{j=0}^{k-1}||\hat{\mathbf{U}}^{j} - \hat{\mathbf{V}}^{j}||\]
holds uniformly on $0\leq t_k \leq T$.
\end{lemma}
\begin{proof}
Let  $R_{1,1}(\tau \Lambda^{\frac{\alpha}{2}}) = \left(4(2I+\tau \Lambda^{\frac{\alpha}{2}})^{-1} -I\right)$. Then  (\ref{eqn:ETDCN}) can be  written as 
\begin{equation*}
\begin{split}
\hat{\mathbf{V}}^{k+1} &= \hat{\mathbf{a}}^{k} + \dfrac{1}{\tau}\Lambda^{-\alpha}\left(R_{1,1}(\tau \Lambda^{\frac{\alpha}{2}}) - I + \tau \Lambda^{\frac{\alpha}{2}}\right)\left[\hat{\mathbf{F}}(\hat{\mathbf{a}}^{k}) - \hat{\mathbf{F}}(\hat{\mathbf{V}}^{k})\right],\\
\hat{\mathbf{a}}^{k} &= R_{1,1}(\tau \Lambda^{\frac{\alpha}{2}})\hat{\mathbf{V}}^{k} +\Lambda^{-\frac{\alpha}{2}}\left( I - R_{1,1}(\tau \Lambda^{\frac{\alpha}{2}})\right)\hat{\mathbf{F}}(\hat{\mathbf{V}}^{k}).
\end{split}
\end{equation*}
Following the same reasoning as in the proof of Lemma (\ref{lem:UnEstimate3}), we obtain
\begin{equation*}
\begin{split}
\hat{\mathbf{a}}^{k-1} &= \left(R_{1,1}(\tau \Lambda^{\frac{\alpha}{2}})\right)^{k}\hat{\mathbf{U}}^{0} + \dfrac{1}{\tau}\Lambda^{-\alpha}(R_{1,1}(\tau \Lambda^{\frac{\alpha}{2}}) - I + \tau \Lambda^{\frac{\alpha}{2}}) \sum_{j=1}^{k-1}\left(R_{1,1}(\tau \Lambda^{\frac{\alpha}{2}})\right)^j\left(\hat{\mathbf{F}}(\hat{\mathbf{a}}^{k-j-1}) - \hat{\mathbf{F}}(\hat{\mathbf{V}}^{k-j-1})\right)\\
&~~ + \Lambda^{-\frac{\alpha}{2}}\left(R_{1,1}(\tau \Lambda^{\frac{\alpha}{2}}) - I\right)\sum_{j=0}^{k-1}\left(R_{1,1}(\tau \Lambda^{\frac{\alpha}{2}})\right)^j\hat{\mathbf{F}}(\hat{\mathbf{V}}^{k-j-1}).
\end{split}
\end{equation*}
Then,
\begin{equation*}
\begin{split}
\hat{\mathbf{b}}^{k-1} - \hat{\mathbf{a}}^{k-1} &= \left(e^{-k\tau  \Lambda^{\frac{\alpha}{2}}} - \left(R_{1,1}(\tau \Lambda^{\frac{\alpha}{2}})\right)^k\right)\hat{\mathbf{U}}^{0} + \dfrac{1}{\tau}\sum_{j=1}^{k-1}\left[\Lambda^{-\alpha}\left(e^{-(j+1)\tau  \Lambda^{\frac{\alpha}{2}}} -\left(R_{1,1}(\tau \Lambda^{\frac{\alpha}{2}})\right)^{j+1}\right)\right.\\
&~~\left.-\Lambda^{-\alpha}\left(e^{-j\tau  \Lambda^{\frac{\alpha}{2}}} - \left(R_{1,1}(\tau \Lambda^{\frac{\alpha}{2}})\right)^j\right)  + \tau  \Lambda^{-\frac{\alpha}{2}}\left(e^{-j\tau  \Lambda^{\frac{\alpha}{2}}} - \left(R_{1,1}(\tau \Lambda^{\frac{\alpha}{2}})\right)^j\right)\right]\\
&~~\times \left(\hat{\mathbf{F}}(\hat{\mathbf{b}}^{k-j-1}) - \hat{\mathbf{F}}(\hat{\mathbf{U}}^{k-j-1})\right)\\
&~~+ \dfrac{1}{\tau}\Lambda^{-\alpha}\left( R_{1,1}(\tau  \Lambda^{\frac{\alpha}{2}}) - I + \tau  \Lambda^{\frac{\alpha}{2}}\right)\sum_{j=1}^{k-1}\left(R_{1,1}(\tau \Lambda^{\frac{\alpha}{2}})\right)^j\left(\hat{\mathbf{F}}(\hat{\mathbf{b}}^{k-j-1}) - \hat{\mathbf{F}}(\hat{\mathbf{a}}^{k-j-1})\right)\\
&~~-\dfrac{1}{\tau}\Lambda^{-\alpha}\left( R_{1,1}(\tau  \Lambda^{\frac{\alpha}{2}}) - I + \tau  \Lambda^{\frac{\alpha}{2}}\right)\sum_{j=1}^{k-1}\left(R_{1,1}(\tau \Lambda^{\frac{\alpha}{2}})\right)^j\left(\hat{\mathbf{F}}(\hat{\mathbf{U}}^{k-j-1}) - \hat{\mathbf{F}}(\hat{\mathbf{V}}^{k-j-1})\right)\\
&~~+ \Lambda^{-\frac{\alpha}{2}} \sum_{j=0}^{k-1}\left(\left(e^{-(j+1)\tau  \Lambda^{\frac{\alpha}{2}}} -\left(R_{1,1}(\tau \Lambda^{\frac{\alpha}{2}})\right)^{j+1}\right) - \left(e^{-j\tau  \Lambda^{\frac{\alpha}{2}}} - \left(R_{1,1}(\tau \Lambda^{\frac{\alpha}{2}})\right)^j\right)  \right)\hat{\mathbf{F}}(\hat{\mathbf{U}}^{k-j-1})\\
&~~+\Lambda^{-\frac{\alpha}{2}}\left( R_{1,1}(\tau  \Lambda^{\frac{\alpha}{2}}) - I \right)\sum_{j=0}^{k-1}\left(R_{1,1}(\tau \Lambda^{\frac{\alpha}{2}})\right)^j\left(\hat{\mathbf{F}}(\hat{\mathbf{U}}^{k-j-1}) - \hat{\mathbf{F}}(\hat{\mathbf{V}}^{k-j-1})\right).
\end{split}
\end{equation*}
Since $\hat{\mathbf{F}}$ is Lipschitz continuous, we have
\begin{equation*}
\begin{split}
||\hat{\mathbf{b}}^{k-1} - \hat{\mathbf{a}}^{k-1}|| &\leq K\tau^2||\hat{\mathbf{U}}^0|| + K\tau^2\sum_{j=1}^{k-1}||\hat{\mathbf{b}}^{k-j-1} - \hat{\mathbf{U}}^{k-j-1}||  + K\tau\sum_{j=1}^{k-1}||\hat{\mathbf{b}}^{k-j-1} - \hat{\mathbf{a}}^{k-j-1}|| \\
&~~+  K\tau\sum_{j=0}^{k-1}||\hat{\mathbf{U}}^{k-j-1} - \hat{\mathbf{V}}^{k-j-1}|| +  K\tau^3\sum_{j=0}^{k-1}||\hat{\mathbf{U}}^{k-j-1}|| + K\tau^3D\\
&\leq  K\tau^2||\hat{\mathbf{U}}^0|| + K\tau^2\sum_{j=0}^{k-2}||\hat{\mathbf{b}}^j - \hat{\mathbf{U}}^{j}|| + K\tau\sum_{j=0}^{k-2}||\hat{\mathbf{b}}^{j} - \hat{\mathbf{a}}^{j}|| \\
&~~+  K\tau\sum_{j=0}^{k-1}||\hat{\mathbf{U}}^{j} - \hat{\mathbf{V}}^{j}|| +  K\tau^3\sum_{j=0}^{k-1}||\hat{\mathbf{U}}^{j}|| + K\tau^3D.
\end{split}
\end{equation*}
Using Lemma \ref{lem:UnEstimate2},
\begin{equation*}
\begin{split}
||\hat{\mathbf{b}}^{k-1} - \hat{\mathbf{a}}^{k-1}|| & \leq  K\tau^2||\hat{\mathbf{U}}^0|| + K\tau^3\sum_{j=0}^{k-1}||\hat{\mathbf{U}}^{j}|| + K\tau\sum_{j=0}^{k-2}||\hat{\mathbf{b}}^{j} - \hat{\mathbf{a}}^{j}|| \\
&~~+  K\tau\sum_{j=0}^{k-1}||\hat{\mathbf{U}}^{j} - \hat{\mathbf{V}}^{j}||  + K\tau^3D
\end{split}
\end{equation*}
and the result follows by Gronwall's lemma for the discrete case\cite[pp. 8]{Emmrich1999}.
\end{proof}
\begin{theorem}
\label{thm:ETDCNEstimate}
Suppose $\hat{\mathbf{F}}(\hat{\mathbf{U}})$ is Lipschitz continuous, then the following error bound 
\[||\hat{\mathbf{U}}^{k} - \hat{\mathbf{V}}^{k}|| \leq K\tau^2\left(||\hat{\mathbf{U}}^{0}|| + \tau D\right)\]
holds uniformly on $0\leq t_k \leq T$.
\end{theorem}
\begin{proof}
We begin by noting that 
\begin{equation*}
\begin{split}
||\hat{\mathbf{U}}^{k} - \hat{\mathbf{V}}^{k}|| & \leq ||\hat{\mathbf{b}}^{k-1}-\hat{\mathbf{a}}^{k-1}|| + \dfrac{1}{\tau}||\Lambda^{-\alpha}\left(e^{-\tau  \Lambda^{\frac{\alpha}{2}}} - R_{1,1}(\tau  \Lambda^{\frac{\alpha}{2}})\right)||\cdot||\hat{\mathbf{F}}(\hat{\mathbf{b}}^{k-1}) - \hat{\mathbf{F}}(\hat{\mathbf{U}}^{k-1})||\\
&~~+\dfrac{1}{\tau}||\Lambda^{-\alpha}\left(R_{1,1}(\tau  \Lambda^{\frac{\alpha}{2}}) -I + \tau  \Lambda^{\frac{\alpha}{2}}\right)||\cdot||\hat{\mathbf{F}}(\hat{\mathbf{b}}^{k-1}) - \hat{\mathbf{F}}(\hat{\mathbf{a}}^{k-1})||\\
&~~+\dfrac{1}{\tau}||\Lambda^{-\alpha}\left(R_{1,1}(\tau  \Lambda^{\frac{\alpha}{2}}) -I + \tau  \Lambda^{\frac{\alpha}{2}}\right)||\cdot||\hat{\mathbf{F}}(\hat{\mathbf{U}}^{k-1}) - \hat{\mathbf{F}}(\hat{\mathbf{V}}^{k-1})||\\
& \leq  K\tau^2||\hat{\mathbf{U}}^{0}|| + K\tau^3\sum_{j=0}^{k-2}||\hat{\mathbf{U}}^{j}|| + K\tau^3D +  K\tau\sum_{j=0}^{k-1}||\hat{\mathbf{U}}^{j} - \hat{\mathbf{V}}^{j}|| \\
&~~+ K\tau^2||\hat{\mathbf{b}}^{k-1} - \hat{\mathbf{U}}^{k-1}|| +  K\tau ||\hat{\mathbf{a}}^{k-1} - \hat{\mathbf{b}}^{k-1}|| + K\tau||\hat{\mathbf{U}}^{k-1} - \hat{\mathbf{V}}^{k-1}||.
\end{split}
\end{equation*}
By Lemmas \ref{lem:UnEstimate2} and \ref{lem:UnEstimate3} and Gronwall's lemma for the discrete case\cite[pp. 8]{Emmrich1999}, the final result is obtained.
\end{proof}
}

\section{Numerical Results}

\blue{The following reaction-diffusion system has been generalised for a fractional derivative.
In the classical case $\alpha = 2$, it reduces to Equation 2 in Garvie's analysis\cite{garvie2007finite}.
}
\begin{align}
\begin{aligned}
\frac{\partial u}{\partial t} &= \Delta^{\alpha/2} u + u(1-u)-vh(au),\\
\frac{\partial v}{\partial t} &= \delta \Delta^{\alpha/2} v + bvh(au)-cv.
\end{aligned}
\label{eqn:predprey}
\end{align}
\blue{The population densities of the prey and predator species are denoted by $u$ and $v$, respectively.} Each population density is considered over time $t$ and vector position $\vec{x}$. The Laplace operator is designated as $\Delta$. The functional response $h$ is defined as the rate of prey consumption per predator relative to the maximum consumption; it increases strictly on $[0,\infty ]$ and satisfies the following conditions:
\begin{align}
\begin{aligned}
h(0)&=0, \\
\lim_{x\to\infty}h(x)&=1.
\end{aligned}
\end{align}
\blue{
The remaining variables $a$, $b$, $c$ and $\delta$ are parameters of the reaction-diffusion system and are strictly positive.
Since the general predator-prey system has been scaled to a non-dimensional form, these parameters do not directly correspond with any intrinsic rates or limits of the biological system.
}
This model does not include an abiotic component and does not account for stochastic factors.
We use the type II functional response proposed by Holling\cite{holling1965functional}:
\begin{equation}
h(\eta)=\frac{\eta}{1+\eta},\text{ }\eta=au.
\end{equation}

\blue{
The following initial and boundary conditions
are associated with this model:
\begin{equation}
\begin{aligned}
\label{eqn:initialBoundaryConditions}
u(x,0) &= f(x), ~~ x\in[a,b]\\
u(0,t) &=0,~~ t \in [0,T]
\end{aligned}
\end{equation}
In our analysis, we assume that positive species densities may exist only within the domain $\Omega$.
Species are not able to leave this `habitat' or grow outside of it, so zero-flux boundary conditions are used.
For a detailed explanation on specific preconditions and implementation,
we refer the reader to Appendix A of \cite{garvie2007finite}.
Additionally, given the model's basis in the natural sciences, conditions that result in negative densities are not considered.
}

\begin{figure}
\centering
\begin{subfigure}{\textwidth}
\includegraphics[width=1.0\textwidth]{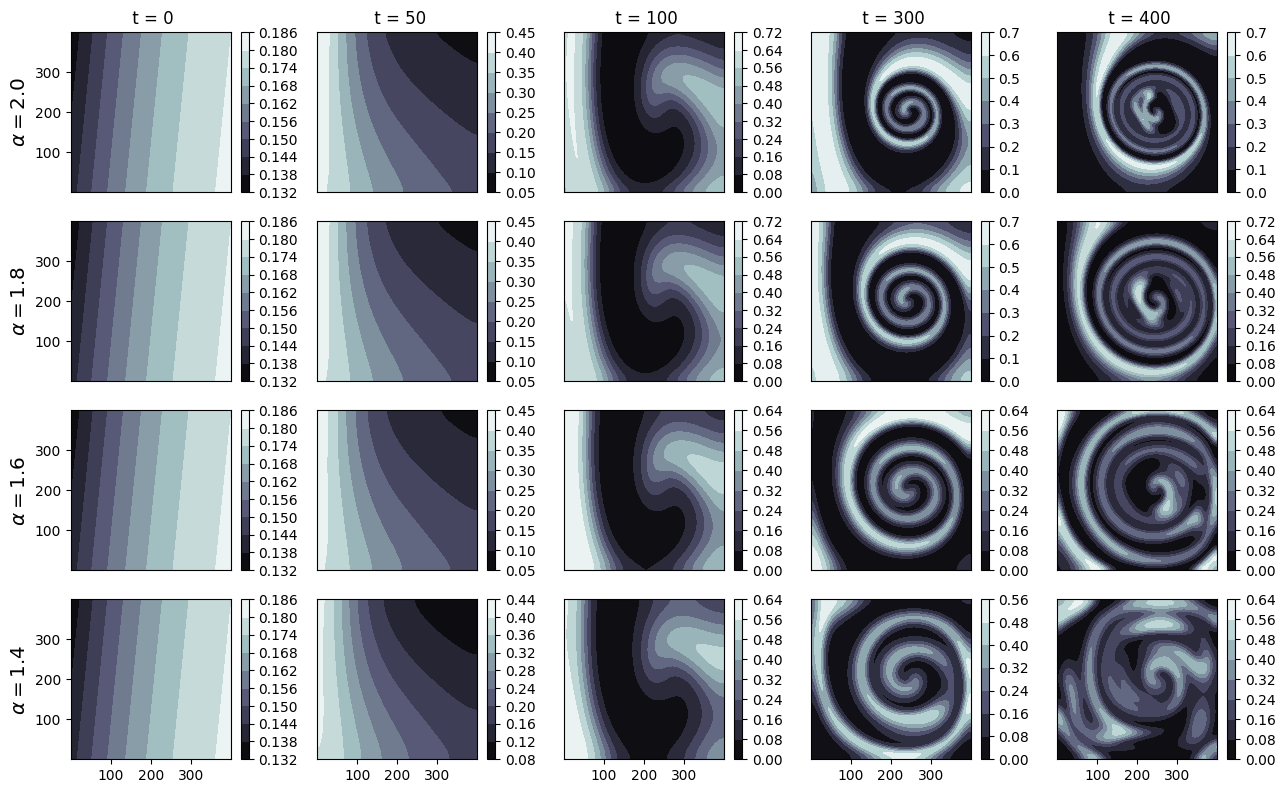}
\caption{Density of U}
\end{subfigure}
\begin{subfigure}{\textwidth}
\includegraphics[width=1.0\textwidth]{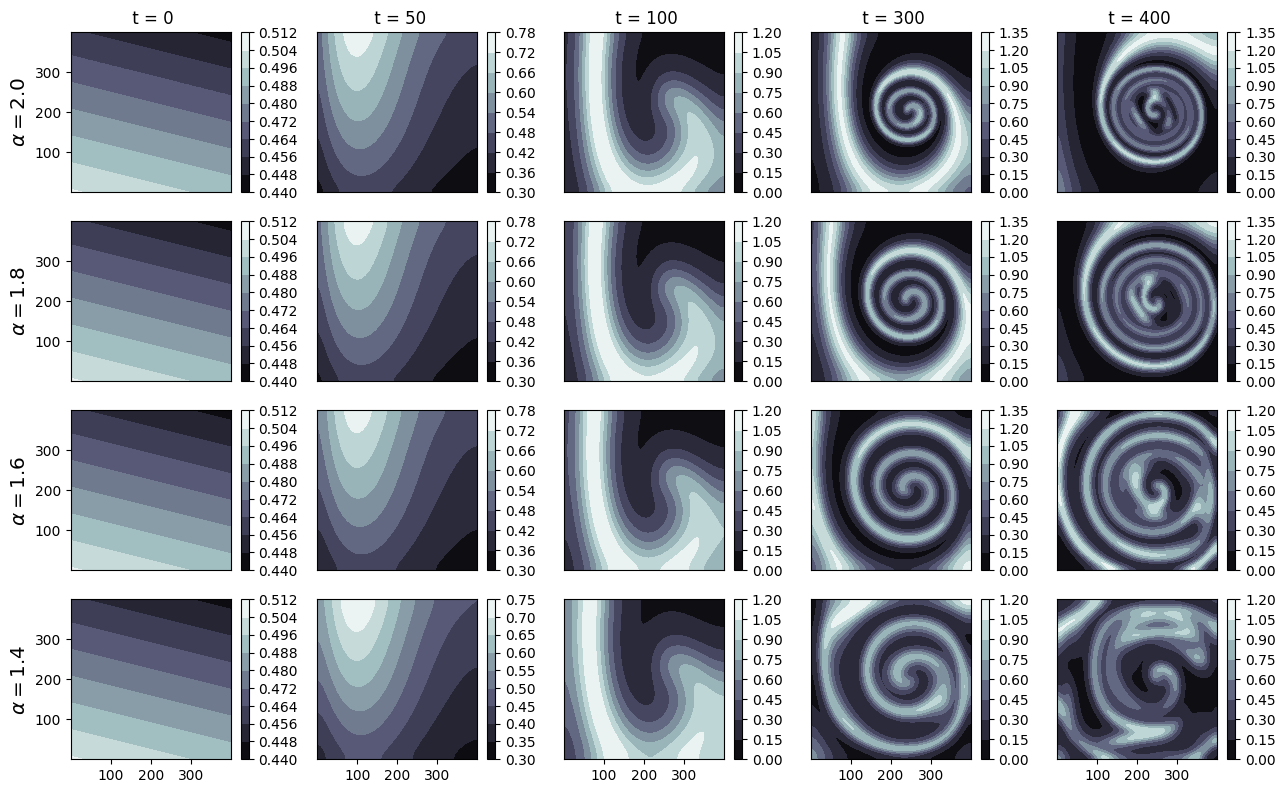}
\caption{Density of V}
\end{subfigure}
\caption{Holling's functional response (Condition A)}
\label{Hollingsa}
\end{figure}

\begin{figure}
\centering
\begin{subfigure}{\textwidth}
\includegraphics[width=1.0\textwidth]{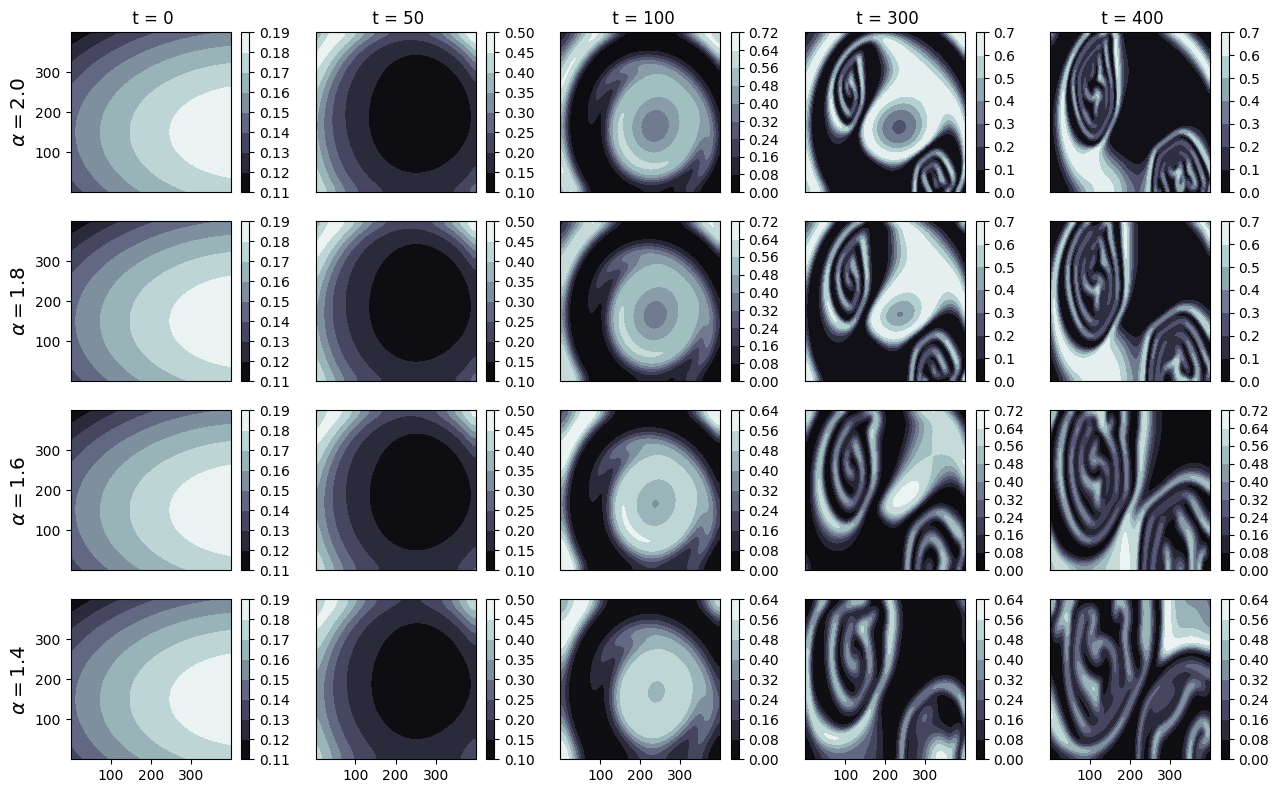}
\caption{Density of U}
\end{subfigure}
\begin{subfigure}{\textwidth}
\includegraphics[width=1.0\textwidth]{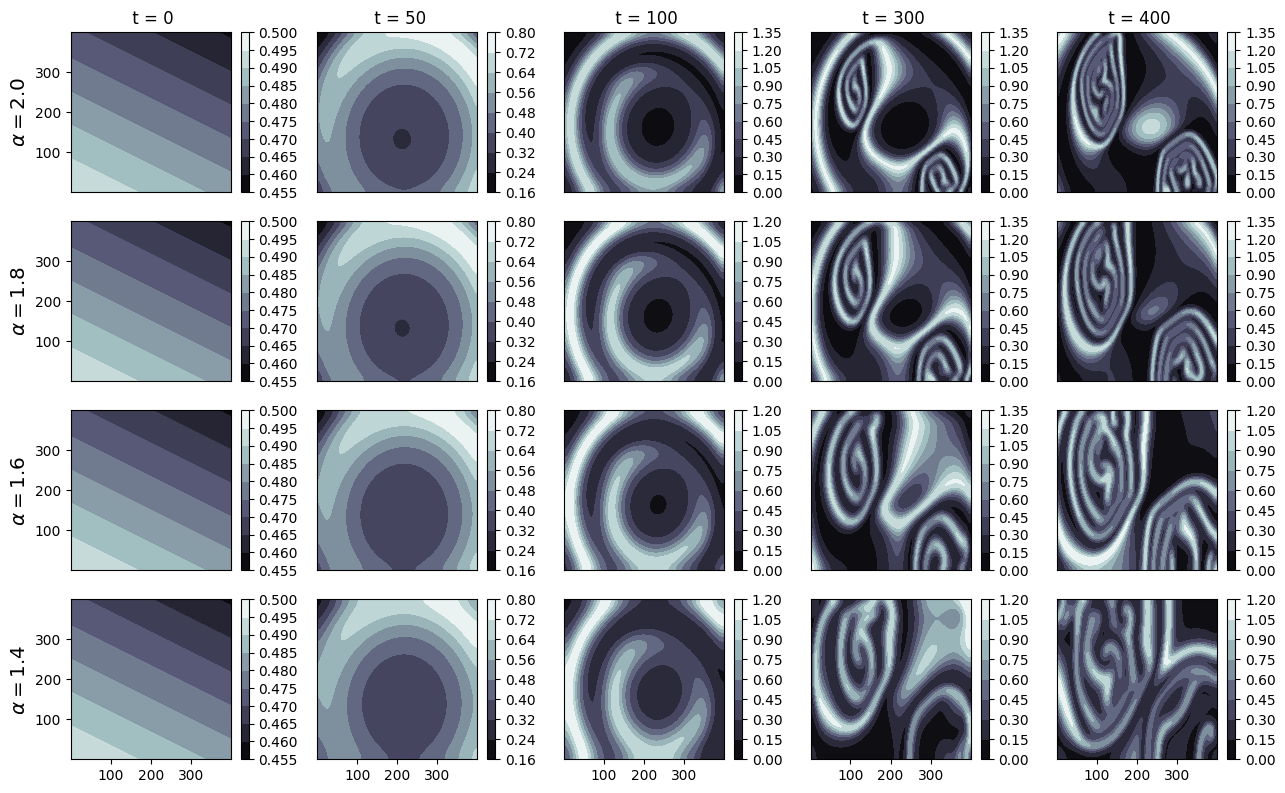}
\caption{Density of V}
\end{subfigure}
\caption{Holling's functional response (Condition B)}
\label{Hollingsb}
\end{figure}

With the proposed methods, we evaluate \autoref{eqn:predprey}
with the parameters
$a=1/0.4$,
$b=2.0$,
\blue{
$c=0.6$,
$\delta = 1$,
}
$u^*= 6/35$, and
$v^*= 116/245$.
\blue{
For these parameters, two different sets of initial conditions are considered.
}
For the set of initial conditions
\begin{align}
U^0_{i,j} &= u^* - 2\times 10^{-7} (x_i - 0.1y_j - 225)(x_i - 0.1y_j - 675),\\
V^0_{i,j} &= v^* - 3\times 10^{-5} (x_i - 450) - 1.2 \times 10^{-4} (y_j - 150),
\end{align}
the results are shown in \autoref{Hollingsa}.
\autoref{Hollingsb} shows the outcome of 
\begin{align}
U^0_{i,j} &= u^* - 2\times 10^{-7} (x_i - 180)(x_i-720) - 6\times10^{-7}(y_j-90)(y_j-210),\\
V^0_{i,j} &= v^* - 3\times 10^{-5} (x_i - 450) - 6 \times 10^{-5} (y_j - 135).
\end{align}

For both sets of initial conditions, we observe that complex, semi-stable patterns emerge from a simple gradient.
As the value of $\alpha$ becomes lower, spreading behavior becomes wider, which is characteristic for fractional-order diffusion\cite{baeumer2007fractional}.
This faster spreading means that wavefronts will approach the boundaries sooner, which alters the pattern development. 

\section{Conclusion}

In this paper, we extend Garvie's work\cite{garvie2007finite} to systems of equations using the a scheme based on the work of Baeumer \textit{et al.}\cite{baeumer2007fractional} and apply the Fourier spectral methods proposed by Bueno-Orovio \textit{et al.}\cite{bueno2014fourier}.
Fractional models serve as a generalisation of existing models and are suitable for anomalous diffusion and nonlocality.
The solution profiles of fractional RDEs have thicker tails than traditional RDEs, so they are better suited \blue{for modelling} species density and biological invasions.

\blue{
\section{Author Contributions}
\textbf{A.~P.~Harris:} Conceptualization, Investigation, Methodology, Software, Visualization, Writing
\textbf{T.~A.~Biala:} Formal analysis, Investigation, Methodology, Software, Visualization, Writing
\textbf{A.~Q.~M.~Khaliq:} Conceptualization, Supervision
}

\appendix

\clearpage
\bibliography{mybib}

\end{document}